\documentclass[12pt]{amsart}
\usepackage{amscd,amssymb,epsf}
\usepackage[mathscr]{eucal}
\newtheorem{thm}{Theorem}[subsection]
\newtheorem*{thm*}{Theorem}
\newtheorem{lemma}[thm]{Lemma}
\newtheorem*{lemma*}{Lemma}

\newtheorem{proposition}[thm]{Proposition}
\newtheorem{corollary}[thm]{Corollary}
\newtheorem*{corollary*}{Corollary}
\newtheorem*{claim*}{Claim}

\newtheorem*{remark*}{Remark}
\numberwithin{equation}{section}

\begin{document}
\newcommand{\R}{{\mathbb R}}
\newcommand{\C}{{\mathbb C}}
\newcommand{\Z}{{\mathbb Z}}
\newcommand{\B}{{\mathbb B}}
\renewcommand{\P}{{\mathbb P}}
\renewcommand{\Re}{{\mathsf{Re}}}
\renewcommand{\Im}{{\mathsf{Im}}}
\newcommand{\Id}{{\mathbb I}}
\newcommand{\Fix}{{\mathsf{Fix}}}
\newcommand{\slt}{{\mathsf{SL}(2,\C)}}			
\newcommand{\sltr}{{\mathsf{SL}(2,\R)}}			
\newcommand{\pslt}{{\mathsf{PSL}(2,\C)}}
\newcommand{\sut}{{\mathsf{SU}(2)}}
\renewcommand{\put}{{\mathsf{PU}(2)}}
\newcommand{\tr}{{\mathsf{tr}}}
\newcommand{\Ad}{{\mathsf{Ad}}}
\newcommand{\dev}{{\mathsf{dev}}}
\newcommand{\gslt}{{\mathfrak{sl}(2,\C)}}		
\newcommand{\gsltr}{{\mathfrak{sl}(2,\R)}}		
\newcommand{\glt}{{\mathfrak{gl}(2,\C)}}		
\newcommand{\ra}{\rangle}
\newcommand{\la}{\langle}
\newcommand{\Hom}{{\mathsf{Hom}}}
\newcommand{\Ham}{{\mathsf{Ham}}}
\newcommand{\hpg}{{\Hom(\pi,\slt)}}
\newcommand{\hpgg}{{\hpg/\hspace{-3pt}/\slt}}
\newcommand{\Ou}{{\mathsf{Out}(\pi)}}
\newcommand{\Au}{{\mathsf{Aut}(\pi)}}
\newcommand{\In}{{\mathsf{Inn}(\pi)}}
\renewcommand{\ss}{{\mathcal{S}}}
\newcommand{\uu}{{\mathcal{U}}}
\newcommand{\tM}{{\tilde{M}}}
\newcommand{\qf}{{\mathcal{QF}(M)}}
\newcommand{\Teich}{{\mathfrak{T}(M)}}
\newcommand{\bTeich}{{\bar{\mathfrak{T}}(M)}}
\newcommand{\cpo}{{\mathbb{CP}^1}}
\newcommand{\cho}{{{\mathsf{H}}_{\C}^1}}
\newcommand{\puoo}{{\mathsf{PU}(1,1)}}
\newcommand{\hol}{{\mathsf{hol}}}
\newcommand{\Mod}{{\mathfrak{M}(M)}}
\newcommand{\gadf}{\gslt_{\Ad \phi}}
\newcommand{\pp}{\mathcal{P}}
\newcommand{\ooo}{\mathcal{O}}
\newcommand{\gp}{\Gamma_{\pp}}
\newcommand{\fc}{{\mathcal F}^{irr}(E)} %
\newcommand{\ga}{{\mathcal G}(E)} 

\title[$\slt$-Character Varieties ]
{The complex-symplectic geometry of $\slt$-characters over surfaces}
\author{William M.~Goldman}
\address{ Mathematics Department,
University of Maryland, College Park, MD  20742 USA  }
\email{ wmg@math.umd.edu }
\thanks{The author gratefully acknowledges partial support from
National Science Foundation grant DMS-0103889.}
\date{\today}

\begin{abstract}
The $\slt$-character variety $X$ of a closed surface $M$ enjoys a
natural complex-symplectic structure invariant under the mapping class
group $\Gamma$ of $M$.
Using the ergodicity
of $\Gamma$ on the $\sut$-character variety, we deduce that every
$\Gamma$-invariant meromorphic function on $X$ is constant.
The trace functions of closed curves on $M$
determine regular functions which generate complex Hamiltonian
flows. For simple closed curves, these complex Hamiltonian flows arise
from holomorphic flows on the representation variety generalizing the
Fenchel-Nielsen twist flows on Teichm\"uller space and the complex
quakebend flows on quasi-Fuchsian space. Closed curves in the complex
trajectories of  these flows lift to paths in the deformation space
$\cpo(M)$ of complex-projective structures between different $\cpo$-structures
with the same holonomy (grafting). 
If $\pp$ is a pants decomposition, then the trace map $\tau_\pp:X
\longrightarrow \C^\pp$ defines a holomorphic completely integrable
system. Furthermore, if $\gp$ is the subgroup of $\Gamma$ preserving
$\pp$, then every $\gp$-invariant holomorphic function
$X\longrightarrow\C$ factors through $\tau_\pp$. This holomorphic
integrable system is related to the complex Fenchel-Nielsen coordinates
on quasi-Fuchsian space $\qf$ developed by Tan and Kourouniotis, and relate
to recent formulas of Platis and Series on complex-length functions
and complex twist flows on $\qf$.
\end{abstract}
\dedicatory{Dedicated to M.S. Raghunathan on his sixtieth birthday}
\maketitle
\tableofcontents
\section*{Introduction}

The {\em $\slt$-character variety\/} of a closed surface $M$ with
fundamental group $\pi$ is the space of equivalence classes of
representations $\phi$ of the fundamental group $\pi$ of $M$ into $\slt$.
The Teichm\"uller space $\Teich$, the moduli space $X_U$ of
irreducible flat unitary $\sut$-bundles, the deformation space
$\cpo(M)$ of $\cpo$-structures on $M$, and the quasi-Fuchsian space
$\qf$ all lie in the smooth stratum $X$ of the $\slt$-character
variety.  These spaces all share several common features: a symplectic
geometry derived from the topology of $M$, and a compatible action of
the mapping class group $\Gamma$ of $M$.  This paper
investigates these structures in terms of the $\Gamma$-invariant
complex-symplectic structure on $X$.

A {\em complex-symplectic structure\/} on a complex manifold
is a nondegenerate closed holomorphic exterior 2-form. The $\slt$-character
variety is an affine variety defined over $\C$ whose set of smooth $\C$-points
is a complex-symplectic manifold of complex dimension $-3\chi(M)$ 
whose complex-symplectic structure is invariant under $\Gamma$. 

In~\cite{FenchelNielsen}, Fenchel and Nielsen develop a set of coordinates
for $\Teich$ based on hyperbolic geometry. This set of coordinates is based
on a {\em pants decomposition\/} 
$\pp = \{\alpha_1,\dots,\alpha_N\}$ on $M$, that is a set of disjoint simple
closed curves $\alpha_i$ cutting $M$ into three-holed spheres (``pants''). 
In a hyperbolic structure on $M$, the curves $\alpha_i$ are represented by 
disjoint simple closed geodesics and their lengths define a function
\begin{equation*}
l_\pp: \Teich \longrightarrow (\R_+)^N 
\end{equation*}
which consitute half of the Fenchel-Nielsen coordinates on $\Teich$.
According to Wolpert~\cite{WP}, these functions Poisson-commute and
define a completely integrable Hamiltonian system. In the language
of classical mechanics, these are the action variables, for which the
Fenchel-Nielsen twist vector fields define angle variables.
In particular $l_\pp$ is the moment map for a Hamiltonian $\R^N$-action.
Choosing a section $\sigma$ to $l_\pp$ determines a symplectomorphim
\begin{align*}
(\R_+)^N \times \R^N  & \longrightarrow \Teich \\
(\lambda; t) & \longmapsto
\xi_t \sigma(\lambda)
\end{align*}
where  $\xi_t$ is the Hamiltonian $\R^N$-action defined by
$l_\pp$. See Wolpert~\cite{FenchelNielsen,Symplectic,WP} for details.

This picture motivated the study of a symplectic geometry of 
moduli spaces $\Hom(\pi,G)/G$ developed in
\cite{Nature,Inv}. The results there extend directly to the
complex-symplectic geometry of deformation spaces $\Hom(\pi,G)/G$
where $G$ is a {\em complex Lie group\/} with an $\Ad$-invariant
complex-orthogonal structure $\B$ on its Lie algebra. In this paper we
consider only the special case $G=\slt$ where $\B$ is the trace form.
Corresponding to an element $\alpha\in\pi$ is the function
$f_\alpha:X\longrightarrow\C$ associating to the equivalence class of
$\phi$ the trace of $\phi(\alpha)$. When $\alpha$ is represented by a
{\em simple\/} closed curve $A$, then the complex-Hamiltonian vector field
$\Ham(f_\alpha)$ generates a flow which is covered by a {\em complex
twist flow\/} on $\hpg$. We study this flow, computing its periods,
and relating it to the action on $X$ of the Dehn twist about $A$.
Following \cite{Nature,Inv}, we relate the complex-symplectic geometry
and the Hamiltonian actions to the Fenchel-Nielsen twist flows on
$\Teich$ and the complex earthquakes and bending deformations on
$\cpo(M)$ and $\qf$. 

In the presence of a conformal structure on $M$, the moduli spaces
$\Hom(\pi,G)/G$ admit stronger structures (a K\"ahler structure when
$G$ is compact, and hyper-K\"ahler when $G$ is complex), but these
stronger structures fail to be $\Gamma$-invariant.  However the
symplectic (and complex-symplectic) structures are $\Gamma$-invariant.
The symplectic structure defines a $\Gamma$-invariant measure.
By \cite{Erg} and Pickrell-Xia~\cite{PickrellXia}, the
resulting measure is ergodic under $\Gamma$ when $G$ is a compact Lie
group.  This has the following consequence for the holomorphic
geometry when $G$ is complex:

\begin{thm*} 
There are no nonconstant $\Gamma$-invariant meromorphic functions on $X$.
\end{thm*}

The proof uses the inclusion of the set $X_U$ of irreducible unitary
characters in $X$, and the ergodicity of the action of $\Gamma$ on $X_U$
(\cite{Erg}). A key idea in the proof is the action of the subgroup $\gp$
preserving a pants decomposition $\pp$ of $M$. The group $\gp$ is a free
abelian group freely generated by the Dehn twists about the curves in $\pp$.
This $\Z^N$-action lies in a Hamiltonian $\R^N$-action. The map which
associates to $[\phi]$ the collection of traces 
\begin{equation*}
X_U \longrightarrow \R^N 
\end{equation*}
is a moment map for the $\R^N$-action and 
is also the ergodic decomposition for the $\Z^N$-action. The holomorphic
analog is:

\begin{thm*} 
Every $\gp$-invariant meromorphic function on $X$ factors through the map
\begin{equation*}
X \longrightarrow \C^N 
\end{equation*}
which associates to a character its values on the curves in $\pp$.
\end{thm*}


The complex twist flows have been extensively studied in the {\em
quasi-Fuchsian space $\qf$,\/} which is the open subset of the
$\slt$-character variety comprising equivalence classes of
quasi-Fuchsian representations. In this case, the complex twist flows
correspond geometrically to quake-bending pleated surfaces in quasi-Fuchsian
hyperbolic 3-manifolds, that is, composing Fenchel-Nielsen twist flows
(earthquakes) with bending deformations (Epstein-Marden~\cite{EpsteinMarden}).
These quakebends are defined more generally for geodesic laminations, 
although we only consider deformations supported on simple closed curves
here. The complex Fenchel-Nielsen coordinates of
Kourouniotis~\cite{Kourouniotis1,Kourouniotis2,Kourouniotis3})
and Tan~\cite{Tan} are holomorphic Darboux coordinates for the
com\-plex-sym\-plec\-tic structure.
We recover results of Platis~\cite{Platis} expressing the symplectic
duality between the complex twist flows and the complex length functions,
and the formula of Series~\cite{Series} for the derivative of a complex
length function under a twist flow in terms of the complex-symplectic
geometry of $X$.

More generally the complex twist flows are defined for 
{\em $\cpo$-structures,\/} that is, geometric structures with coordinates
modelled on $\cpo$ with coordinate changes in $\pslt$. Let $\cpo(M)$
denote the deformation space of $\cpo$-structures on $M$. 
Using the local biholomorphicity of the
holonomy mapping
\begin{equation*}
\hol:  \cpo(M) \longrightarrow \hpgg.
\end{equation*}
one obtains complex-symplectic structures, complex length functions,
and Hamiltonian complex twist flows on $\cpo(M)$. The complex twist
flows on $\cpo$-structures can be described geometrically by inserting
annuli into a $\cpo$-manifold split along a simple closed curve which
is locally circular. This is a special case of the {\em grafting\/}
construction considered in Tanigawa~\cite{Tanigawa} and
McMullen~\cite{McMullen}. In particular closed curves in the complex
trajectory of a complex twist flow lift to paths between different
$\cpo$-structures with the same holonomy (Maskit and Hejhal, see
Goldman~\cite{Fuchsian}). The holomorphic properties of the grafting
construction are discussed in McMullen~\cite{McMullen},
Tanigawa~\cite{Tanigawa} and Scannell-Wolf~\cite{ScannellWolf}.

\section{Representation varieties and character varieties}

Let $M$ be a closed oriented surface with fundamental group $\pi$.
Let $\hpg$ denote the complex affine variety consisting of
homomorphisms $\pi\longrightarrow\slt$. (It is irreducible, by
Goldman~\cite{TopComps}; see 
Benyash-Krivets -Chernousov -Rapinchuk~\cite{BCR} for
stronger more general results, and also Li~\cite{JunLi}.)  

The group $\slt$ acts by conjugation on the affine algebraic set $\hpg$.  
Let $\hpg/\slt$ denote the set of $\slt$-orbits of $\hpg$.
Denote the orbit of a representation $\phi$ by
\begin{equation*}
[\phi]\in \hpg.
\end{equation*}

\subsection{Stable and semistable points}
A homomorphism 
\begin{equation*}
\phi\in\hpg 
\end{equation*}
is {\em irreducible\/} if it leaves invariant
no proper linear subspace of $\C^2$. Equivalently, $\phi$ is irreducible
if the corresponding projective action fixes no point in $\cpo$.
Irreducible homomorphisms
are the {\em stable points\/} of $\hpg$ (with respect to the $\slt$-action).
The subset of $\hpg^s \subset  \hpg$ consisting of irreducible homomorphisms 
is Zariski-open and nonsingular.  $\hpg^s$ is a complex manifold
of complex dimension $6g-3$ (see~\cite{Nature}).

More generally, a homomorphism is {\em reductive\/} if every invariant
subspace possesses an invariant complement. Equivalently, $\phi$ is
reductive if it is either reducible or fixes a pair of distinct points
on $\cpo$ or is {\em central.\/} A central homomorphism maps into the
center $\{\pm\Id\}$ of $\slt$ and acts trivially on $\cpo$. Letting
$\Fix(\phi)$ denote the subset of $\cpo$ fixed by $\phi$, a
homomorphism $\phi$ is reductive if and only if $\Fix(\phi)$ equals
$\emptyset$, a pair of distinct points, or all of $\cpo$.  Reductive
homomorphisms are the {\em semistable points\/} of $\hpg$ (with
respect to the $\slt$-action),comprising the subset 
\begin{equation*}
\hpg^{ss} \subset \hpg.  
\end{equation*}
A {\em semistable\/} (respectively {\em stable\/})
orbit is the orbit of a semistable (respectively stable) point.

Let $\hpgg$ denote the set of $\C$-points of the {\em categorical quotient\/}
of the $\slt$-action on $\hpg$, that is, the variety whose
coordinate ring is the ring of invariants of $\slt$ acting on the
coordinate ring $\C[\hpg]$ of $\hpg$. The natural map
\begin{equation*}
\hpg/\slt \longrightarrow \hpgg 
\end{equation*}
is surjective, but not injective.
However $\hpgg$ may be identified with
the set of semistable orbits.
Thus the inclusion
\begin{equation*}
\hpg^{ss} \hookrightarrow \hpg
\end{equation*}
induces a bijection
\begin{equation*}
\hpg^{ss}/\slt \longleftrightarrow \hpgg.
\end{equation*}
The group $\slt$ acts freely and properly on $\hpg^s$. The quotient 
$X := \hpg^s/\slt$ is thus a $(6g-6)$-dimensional complex manifold which embeds
as a Zariski open subset of the categorical quotient $\hpgg$. Thus the map
\begin{equation*}
X  = \hpg^s/\slt \longrightarrow \hpgg
\end{equation*}
induced by $\hpg^s\hookrightarrow\hpg$ is an embedding onto an open subset.
Thus $X$ is a 
smooth irreducible complex quasi-affine variety which is dense in the
quotient $\hpgg$.

\subsection
{Symplectic geometry of deformation spaces}
By the general construction of \cite{Nature}, $X$ has a natural {\em
complex-symplectic structure $\Omega$,\/} which is $\Gamma$-invariant.
Furthermore $\Gamma$ is algebraic in the sense that there exists an
algebraic tensor field on $\hpg$ inducing $\Omega$. 
(See \cite{Nature} for an explicit formula.)

The Zariski tangent space $T_{\phi}\hpg$ is the space 
\begin{equation*}
Z^1(\pi,\gadf) 
\end{equation*}
of {\em 1-cocycles \/} $\pi\longrightarrow\gadf$ and the tangent space
$T_{\phi}(G\cdot\phi)$ of the $G$-orbit equals the subspace
\begin{equation*}
B^1(\pi,\gadf)\subset Z^1(\pi,\gadf)
\end{equation*}
of {\em 1-coboundaries.\/} (Compare Raghunathan~\cite{Raghunathan}.)
The quotient vector space is the cohomology 
\begin{equation*}
H^1(\pi,\gadf) 
\end{equation*}
which, under the de Rham isomorphism is isomorphic to the cohomology
$H^1(M;V_\phi)$ where $V_\phi$ denotes the flat vector bundle over $M$
corresponding to the $\pi$-module $\gadf$.

Let $\B:\gslt\times\gslt\longrightarrow\C$  be the {\em trace form\/} of 
the standard representation on $\C^2$: if $\alpha,\beta\in\gslt$, then the
inner product is defined as 
\begin{equation*}
\B(\alpha,\beta) := \tr(\alpha\beta).
\end{equation*}
Since $\B$ is $\Ad$-invariant, it defines a bilinear pairing of
$\pi$-modules
\begin{equation*}
\gadf \times \gadf \longrightarrow \C
\end{equation*}
or equivalently flat vector bundles 
\begin{equation*}
V_\phi \times V_\phi \longrightarrow \C.
\end{equation*}
Cup-product defines a bilinear pairing
\begin{equation*}
\Omega_\phi: H^1(M;V_\phi) \times  H^1(M;V_\phi) \longrightarrow \C
\end{equation*}
with coefficients paired by $\B$. Symmetry of $\B$ implies
that $\Omega_\phi$ is skew-symmetric.
Since $\B$ is nondegenerate, $\Omega_\phi$ is nondegenerate.
It follows from \cite{Nature} that $\Omega_\phi$ can be 
expressed as an algebraic tensor on $\hpg$, and thus is a holomorphic
exterior 2-form. By arguments of \cite{Nature} or 
Karshon~\cite{Karshon}, Weinstein~\cite{Weinstein}, 
Guruprasad-Huebschmann-Jeffrey-Weinstein~\cite{GHJW},
$\Omega$ is closed. Thus $\Omega$ is a nondegenerate closed 
holomorphic $(2,0)$-form, that is a {\em complex-symplectic structure.\/}

\subsection
{The mapping class group}
The automorphism group $\Au$ of $\pi$ acts algebraically on $\hpg$
commuting with the action of $\slt$. The normal subgroup $\In$ of
inner automorphisms acts trivially on the quotient $\hpgg$. The
quotient $\Ou = \Au/\In$, acts on $\hpgg$, leaving invariant the
subset $X$.  
Furthermore the mapping class group
$\Gamma := \pi_0(\mathsf{Diff}^+(M))$  of $M$
is isomorphic to $\Ou$ 
by Nielsen~\cite{Nielsen}.  The algebraic complex-symplectic structure
on $X$ is $\Gamma$-invariant.

\subsection
{Ergodicity and its holomorphic analog}
The subset $X_U$ of unitary characters is invariant under $\Gamma$.
Furthermore the complex-symplectic structure restricts to a (real)
symplectic structure on $X_U$ which is the K\"ahler form for a K\"ahler
structure on $X_U$. In particular the Lebesgue measure class on $X_U$
is invariant under $\Gamma$. The main result of \cite{Erg} is
that this action is {\em ergodic,\/} that is every $\Gamma$-invariant 
measurable function on $X_U$ is constant almost everywhere.

Ergodicity no longer holds for the $\slt$-character variety $X$ (see
\S\ref{sec:quasifuchsian}).  However ergodicity on $X_U$ does imply the
following property of holomorphic functions on $X$:

\begin{thm} A $\Gamma$-invariant meromorphic function 
$X\stackrel{h}\longrightarrow\cpo$ is constant.
\end{thm}
\begin{proof}

The restriction of $h$ to $X_U$ is a $\Gamma$-invariant measurable
function, and by the main result of Goldman~\cite{Erg}, $h$ must
be constant almost everywhere. Since $h$ is continuous, it is constant.

Now we argue that in local holomorphic coordinates, $X_U$ is equivalent to 
$\R^n\subset \C^n$ and a holomorphic function constant on $X_U$ must
be globally constant on $X$. (Compare, for example, Lemma 1 of \S 2.3 of 
Platis~\cite{Platis}).
For the reader's convenience, we supply a brief proof.
Let $\uu$ be a nonempty open coordinate neighborhood with 
local holomorphic coordinates $z = (z^1,\dots,z^n)$.
In local holomorphic coordinates, $X_U\cap\uu$ 
is described by 
\begin{equation*}
z\in \R^n\subset \C^n 
\end{equation*}
and $h$ is given by a power series 
\begin{equation*}
h(z) = \sum_{k=0}^\infty a_k z^k 
\end{equation*}
which converges in the nonempty open set $z(\uu)\in\C^n$. 
Since the restriction of $h$ to a $z(\uu)\cap \R^n$ is constant,
$a_k = 0$ for $k>0$ and thus $h$ is constant on $\uu$.
Since $X$ is connected, analytic continuation implies that $h$ must be 
constant. 
\end{proof}
We do not know whether $X-X_U$ admits $\Gamma$-invariant nonconstant 
meromorphic functions.

\section
{The Hamiltonian vector field of a character function}

Corresponding to free homotopy classes $\alpha$ of closed curves on $M$ are
complex regular functions $f_\alpha:X\longrightarrow\C$. 
The complex-symplectic structure associates to these functions
complex Hamiltonian vector fields  $\Ham(f_\alpha)$, which generate
holomorphic local flows on $X$. 
When $\alpha$ corresponds to a simple closed curve, we define
{\em complex twist flows\/} on $\hpg$ which cover these holomorphic local
flows on $X$. We begin with an preliminary section on the traces 
in $\slt$.

\subsection
{The variation of the trace function on $\slt$}
Let 
\begin{equation*}
f:\slt\longrightarrow\C 
\end{equation*}
be the trace function $f(P)=\tr(P)$ and $\B$ the trace form
\begin{align*}
\gslt \times  \gslt &\longrightarrow \C \\
(X,Y) & \longmapsto \tr(XY).
\end{align*}
As in Goldman~\cite{Inv} the differential of $f$ and the
orthogonal structure $\B$ determines a variation function
\begin{equation*}
F:\slt\longrightarrow\gslt  
\end{equation*}
characterized by the identity
\begin{equation*}
\frac{d}{dt} \bigg| _{t=0} f\big(P \exp(tX)\big) = \B\big(F(P),X\big).
\end{equation*}
This function is defined by:
\begin{equation*}
F(P) = P - \frac{\tr(P)}2 \Id 
\end{equation*}
and corresponds to the composition of the inclusion
\begin{equation*}
\slt \hookrightarrow \glt
\end{equation*}
with orthogonal projection 
\begin{equation*}
\glt \hookrightarrow \gslt
\end{equation*}
(orthogonal with respect to the trace form $\B$ on $\gslt$).
Invariance of the trace
\begin{equation*}
f(Q P Q^{-1}) =  f(P)
\end{equation*}
and $\Ad$-invariance of the orthogonal structure $\B$
\begin{equation*}
\B(\Ad(Q)X, \Ad(Q)Y) = \B(X,Y) 
\end{equation*}
implies $\Ad$-equivariance of its variation:
\begin{equation}\label{eq:equivariance}
F(Q P Q^{-1}) =  \Ad(Q)F(P).
\end{equation}
Taking $Q=P$ shows that $F(P)$ lies in the centralizer of $P$ in $\gslt$.
In particular the complex one-parameter subgroup of $\slt$
\begin{equation}\label{eq:zetat}
\zeta_t := \exp\big(tF(P)\big)  
\end{equation}
centralizes $A$ in $\slt$.

For example, if $P$ is the diagonal matrix
\begin{equation*}
\bmatrix \lambda & 0 \\ 0 & \lambda^{-1} \endbmatrix
\end{equation*}
then $f(P) = \lambda + \lambda^{-1}$ and 
\begin{equation*}
F(P) = \frac{\lambda-\lambda^{-1}}2
\bmatrix 1 & 0 \\ 0 & -1 \endbmatrix
\end{equation*}
with corresponding complex one-parameter subgroup:
\begin{equation}\label{eq:diagC}
\zeta_t =  \bmatrix   
e^{t(\lambda-\lambda^{-1})/2} & 0 \\ 
0 & e^{-t(\lambda-\lambda^{-1})/2}
\endbmatrix.
\end{equation}
If $P$ is $\pm$ the unipotent matrix
\begin{equation*}
\bmatrix 1 & 1 \\ 0 & 1 \endbmatrix
\end{equation*}
then $f(P) = \pm 2$ and
\begin{equation*}
F(P) = \bmatrix 0 & 1 \\ 0 & 0 \endbmatrix
\end{equation*}
with corresponding complex one-parameter subgroup:
\begin{equation}\label{eq:unipC}
\zeta_t =  \bmatrix  1 & t \\ 0 & 1 \endbmatrix\end{equation}

The following lemma, whose proof is immediate from the above calculations,
will be needed in the sequel. Recall that for any subset $S\subset\slt$, 
$\Fix(S)$ denotes the subset of $\cpo$ fixed by $S$.
\begin{lemma}\label{lem:fix}
If $P\in\slt$, then $\Fix(\zeta_t) = \Fix(P)$
or $\zeta_t = \pm \Id$
\end{lemma}

\subsection
{Character functions}
Let $\alpha\in\pi$. The {\em character function\/}
\begin{align*}
\hpg & \stackrel{\tilde{f}_\alpha}\longrightarrow \C \\
\phi \qquad & \longmapsto  f(\phi(\alpha))
\end{align*}
is a regular function. Since $f$ is a class function,
$\tilde{f}_\alpha$ is invariant under $\slt$ and defines a regular function
\begin{equation*}
f_\alpha: X \longrightarrow \C.
\end{equation*}
Let $[\phi]\in X$. Then the
tangent space to $X$ equals $H^1(M;V_\phi)$, so the 
Hamiltonian vector field $\Ham(f_\alpha)$ associates to $[\phi]\in X$
a tangent vector in $H^1(M;V_\phi)$.  

Cap product with the fundamental homology
class $[M]\in H_2(M;\C)$ defines the {\em Poincar\'e duality isomorphism:\/}
\begin{equation*}
\cap [M]: H^1(M;V_\phi) \stackrel{\cong}\longrightarrow H_1(M;V_\phi).
\end{equation*}
Choose a basepoint $x_0\in M$, an isomorphism of the fiber of $V_\phi$
over $x_0$ with $\C^2$, and a representative holonomy homomorphism
\begin{equation*}
\phi_0: \pi_1(M;x_0)  \longrightarrow \slt.
\end{equation*}
Let $s_0\in S^1$ be a basepoint.
Let $\alpha_0:(S^1,s_0)\longrightarrow (M,x_0)$ be a based loop in $M$
corresponding to $\alpha$.  Let $\sigma$ be the parallel section of
the flat vector bundle $\alpha_0^*V_{\phi}$ over $S^1$ which equals
$F(\phi_0(\alpha_0))$ at $s_0$. Then $\sigma$ defines a $V_{\phi}$-valued
1-cycle in $M$, with homology class
\begin{equation*}
[\sigma] \in H_1(M;V_{\phi}).
\end{equation*}
\begin{lemma} The value of the Hamiltonian vector field 
$\Ham(f_\alpha)$ at a point $[\phi]\in X$ is the vector
in 
\begin{equation*}
T_{[\phi]}X \cong H^1(M;V_{\phi}) 
\end{equation*}
corresponding to the Poincar\'e dual of the homology class of the
$V_\phi$-valued cycle $\sigma$:
\begin{align*}
\cap [M]: H^1(M;V_\phi) & \longrightarrow H_1(M;V_\phi) \\
\Ham(f_\alpha)  & \longmapsto [\sigma]. 
\end{align*}
\end{lemma}
Using the above formula and the duality between cup-product and intersections
of cycles, we obtain a formula for the Poisson bracket of trace functions
(\cite{Inv}):
\begin{proposition}\label{prop:PoissonBracket}
Let $\alpha,\beta$ be oriented closed curves meeting transversely
in double points $p_1,\dots,p_k$. For each $p_i$, choose representatives
\begin{equation*}
\phi_i:\pi_1(M;p_i)\longrightarrow \slt.  
\end{equation*}
Let $\alpha_i$ and $\beta_i$
be the elements of $\pi_1(M;p_i)$ representing $\alpha,\beta$ respectively.
Then the Poisson bracket of functions $f_\alpha,f_\beta$ is:
\begin{align*}
\big\{f_\alpha,f_\beta\big\} & =
\Omega\big( \Ham(f_\alpha),\Ham(f_\beta) \big) \\ & =
\sum_{i=1}^k\  \epsilon(p_i;\alpha,\beta) \ 
\B\big(F(\phi_i(\alpha_i)),F(\phi_i(\beta_i))\big)
\end{align*}
where $\epsilon(p_i;\alpha,\beta)$ denotes oriented intersection number.
\end{proposition}
\begin{corollary}\label{cor:disjoint}
If $\alpha,\beta$ are disjoint, then $f_\alpha$ and $f_\beta$ Poisson-commute.
\end{corollary}
By arguments in \cite{Inv} (based on a suggestion of S.\ Wolpert), the
converse holds if one of $\alpha$ or $\beta$ is simple. A purely
topological proof of this fact has recently been given by
Chas~\cite{Chas}.

\subsection
{Complex twist flows}
Suppose that $\alpha$ is represented by a simple closed curve $A$. 
Then a holomorphic $\C$-action on the representation variety $\hpg$ covers
the complex Hamiltonian flow on $X$. 
The {\em complex twist flow\/} $\tau_\alpha$ is the holomorphic
action of $\C$ of $\hpg$ defined as follows. 

There are two cases, depending on whether $A$ separates $M$ or not.
Denote by $M|A$ the compact surface with boundary whose interior is 
homeomorphic to the complement $M-A$. Denote the two components of 
$\partial(M|A)$ by $A_+$ and $A_-$. 
The quotient map $q:M|A \longrightarrow A$ 
results from identifying these components by a homeomorphism
$\eta:A_+ \longrightarrow A_-$. 

\subsubsection
{The nonseparating case.}
If $A$ is nonseparating, then $M|A$ is connected, and choosing a basepoint
in $A_+$ we express 
$\pi$ as an HNN-extension
\begin{equation*}
\pi_1(M|A)\, \star_\iota 
\end{equation*}
as follows.
Choose a basepoint $x_0\in A$, and lift it to a 
basepoint $\tilde{x}_0\in A_+\subset M|A$. Let $\alpha_+$
denote the element of $\pi_1(M|A;\tilde{x}_0)$ corresponding
to $A_+$. Let $\tilde\beta$ denote a simple arc joining
$\tilde{x}_0$ to $\eta(\tilde{x}_0)$ in $M|A$, and
$\alpha_-$ the based loop 
${\tilde\beta}^{-1}\cdot \alpha_+\cdot\tilde\beta$
at the basepoint $\tilde{x}_0$. Then $\eta$ corresponds
to the isomorphism 
\begin{equation*}
\iota:\la\alpha_+\ra\longrightarrow\la\alpha_-\ra
\end{equation*}
between the two subgroups of $\pi_1(M|A)$.
Let $N$ be the normal closure of the element
\begin{equation*}
\beta\alpha_+\beta^{-1}\big(\alpha_-\big)^{-1}
\end{equation*}
of $\pi_1(M|A) \star \la \beta\ra$.
Then the HNN-extension $\pi_1(M|A) \star_\iota$
is defined as the quotient
\begin{equation*}
\pi \cong
\bigg( \pi_1(M|A) \star \la \beta\ra\bigg) / N. 
\end{equation*}
As in \cite{Inv}, for any $\phi\in\hpg$, 
and $t\in\C$, define
\begin{equation}
\xi_t(\phi): \gamma \longmapsto\label{eq:nonsep}
 \begin{cases}
\phi(\gamma) & \text{~if~$\gamma\in\pi_1(M|A)$}  \\
\phi(\beta)\ \zeta_t 
& \text{~if~$\gamma=\beta$.}\end{cases}
\end{equation}
where
\begin{equation*}
\zeta_t = \exp \big(tF(\phi(\alpha_+))\big)\end{equation*}
is the one-parameter subgroup corresponding to $\phi(\alpha_+)$.

\subsubsection
{The separating case.}
If $A$ separates $M$, then let $M_+$ and $M_-$ denote the two
components of $M|A$. Then $\pi$ is the free product of subgroups corresponding
to $\pi_+(M_+)$ and $\pi_1(M_-)$ respectively,  amalgamated over the
images $\la \alpha \ra \longrightarrow \pi_1(M_\pm)$.
Then as in \cite{Inv}, for any $\phi\in\hpg$, $t\in\C$, define
\begin{equation}\label{eq:sep}
\xi^{\alpha}_t(\phi): \gamma \longmapsto \begin{cases}
\phi(\gamma) 
& \text{~if~$\gamma\in\pi_1(M_+)$} \notag\\
\zeta_t \ 
\phi(\gamma)\ 
\zeta_{-t}
& \text{~if~$\gamma\in\pi_1(M_-)$} \end{cases}
\end{equation}
where
\begin{equation*}
\zeta_t = \exp\big(tF(\phi(\alpha)\big)  
\end{equation*}
is the one-parameter subgroup 
corresponding to $\phi(\alpha)$.

\subsubsection
{Irreducibility}
To show that the Hamiltonian flows act on $X$, we need to show that
the complex twist flows preserve the set of irreducible representations.
\begin{lemma}
Suppose that $\phi\in\hpg^s$. Then 
\begin{equation*}
\xi^{\alpha}_t(\phi)\in\hpg^s. 
\end{equation*}
\end{lemma}
\begin{proof}
Let $M'$ be a component of $M|A$. As above, choose basepoints and
identifications to express $\pi$ as an amalgamated free product
or HNN construction with $\pi_1(M')\hookrightarrow \pi_1(M)$. Then 
\begin{equation*}
\xi^{\alpha}_t\phi(\pi) \supset \phi(\pi_1(M')) 
\end{equation*}
implies that if $\phi_{\pi_1(M')}$ is irreducible, then so is
$\xi^{\alpha}_t\phi$. Thus we may assume that 
$\phi_{\pi_1(M')}$ is reducible but $\phi$ is irreducible. 

Suppose first that $A$ is nonseparating. 
We assume $\xi^{\alpha}_t\phi$ is reducible and derive a contradiction.
Let $\Phi = \Fix(\phi(\pi_1(M|A)))$.
Since 
\begin{equation*}
\phi(\alpha_+)\in \phi\big(\pi_1(M|A)\big),  
\end{equation*}
the element $\phi(\alpha_+)$
fixes $\Phi$, and Lemma~\ref{lem:fix} implies that
$\zeta_t$ fixes $\Phi$. Since $\xi^{\alpha}_t\phi(\pi)$ is reducible and is
generated by $\phi(\pi_1(M|A))$ and $\phi(\beta)$, it follows that
\begin{equation*}
\xi^{\alpha}_t\phi(\beta) = \phi(\beta)\zeta_t
\end{equation*}
fixes no element of $\Phi$. Since $\zeta_t$ fixes $\Phi$, this implies
that $\phi(\beta)$ fixes $\Phi$, a contradiction.

Suppose finally that $A$ separates $M$ into two components $M_\pm$. By
the remark above, we may assume that $\phi_{\pi_1(M_\pm)}$ is reducible.
Let 
\begin{equation*}
\Phi_\pm = \Fix(\phi(\pi_1(M_\pm))). 
\end{equation*}
We may assume that
each $\Phi_\pm$ is nonempty but $\Phi_+\cap\Phi_-=\emptyset$.
Now $\phi(\alpha)$ fixes each $\Pi_\pm$, and Lemma~\ref{lem:fix} implies
that $\zeta_t$ fixes each $\Pi_\pm$. 
Since $\phi(\pi_1(M_+))$ and $\zeta_t \phi(\pi_1(M_-))\zeta_t^{-1}$ 
generate $\xi^{\alpha}_t\phi(\pi)$,
\begin{align*}
\Fix(\xi^{\alpha}_t\phi(\pi)) & 
=  \Fix\big(\phi(\pi_1(M_+))\big) \cap \Fix\big(
\zeta_t \phi(\pi_1(M_-))\zeta_t^{-1}\big) \\ &
=  \Fix\big(\phi(\pi_1(M_+))\big) \cap \zeta_t\Fix\big(\phi(\pi_1(M_-)))\big) 
\\ &  = \Phi_+ \cap \zeta_t \Phi_- 
\\ &  = \Phi_+ \cap  \Phi_- = \emptyset
\end{align*}
as desired.
\end{proof}

\subsection{Gauge-Theoretic Interpretation}
These twist flows admit an interpretation in terms of flat
connections.  The character variety $X$ identifies with the quotient
of the space $\fc$ of irreducible flat $\slt$-connections on a
principal $\slt$-bundle $E$ over $M$ by the group $\ga$ of gauge
transformations of $E$.  (Since a principal $\slt$-bundle over $M$
which admits a flat connection is necessarily trivial, we may assume
$E$ is the product bundle.)

Let $A\subset M$ be a simple closed curve, and choose a basepoint  $a_0\in A$.
Choose an orientation on $A$ and let $\alpha$ be a based loop on $M$
corresponding to the orientation on $A$.

Pull $E$ back to a principal $\slt$-bundle $q^*E$ over
$M|A$ by the quotient map $q:M\longrightarrow M|A$. Choose a point
$e_0$ in the fiber of $E$ over $a_0$. Let $a_\pm$ be the two elements
of $q^{-1}(a_0)$ in $A_\pm$ respectively and $e_\pm$ the 
elements in the fibers of $q^*E$ over $a_\pm$ corresponding to $e_0$.

Define the twist flow $\tilde{\xi}_t$ on $\fc$ as follows. 
Let $\nabla\in\fc$ be a flat connection on $E$. Parallel transport
of $e_0$ along $\alpha$ with respect to $\alpha$ defines a holonomy
transformation $\phi_0(\alpha_0)$ as above. Let $\zeta_t$ denote the
corresponding one-parameter subgroup of $\slt$.
There is a one-parameter family of gauge-transformations $g_t$ of $q^*E$
supported in a collar neighborhood $N$ of $A_+\subset M|A$ assuming the
``value'' $\zeta_t$ on $A_+$.

Explicitly, choose a smooth embedding
\begin{equation*}
\psi:[0,1]\times S^1\longrightarrow N  \hookrightarrow M|A
\end{equation*}
mapping $\{1\}\times S^1$ diffeomorphically to $A_+$. Let
$r:[0,1]\longrightarrow S^1$ denote the quotient mapping identifying
$0,1\in [0,1]$. Parallel transport of $e_+$ defines a trivialization
of the pullback of $E$ to $[0,1]\times [0,1]$ so that $\psi^*E$ identifies
with the quotient of 
\begin{equation*}
[0,1]\times [0,1] \times \slt 
\end{equation*}
by the identification
\begin{equation*}
(s,0,h) \longleftrightarrow  (s,1,\phi_0(\alpha_0)h).
\end{equation*}
Define $g_t$ as the gauge transformation which is the identity map on the
complement of $N$ in $M|A$, and equals 
\begin{equation*}
(s,\theta,h) \longrightarrow (s,\theta,\zeta_{s,t}h)
\end{equation*}
in this trivialization.

Since $\zeta_{t}$ centralizes the holonomy of $q^*\nabla$ along $A_+$,
the identification map $\eta:A_-\longrightarrow A_+$ 
identifies the restrictions of $(g_t)^*(q^*\nabla)$ to $A_\pm$.
Thus a unique flat connection $\tilde{\xi}_t(\nabla)$ exists, 
satisfying
\begin{equation*}
q^*\big(\tilde{\xi}_t(\nabla)\big) = 
(g_t)^*(q^*\nabla).
\end{equation*}
This is the orbit of the twist flow on flat connections. Clearly the
gauge transformation $g_t$ does not arise from a gauge transformation
of $E$, unless the holonomy along $A$ is $\pm\Id$. The orbit covers
the orbit of the holonomy of $\nabla$ under the twist flow in
$X \cong \fc/\ga$.

The flow $\tilde\xi$ on $\fc$ depends only on the choice of 
the collar neighborhood $N\subset M|A$ and 
$\psi:[0,1]\times S^1\longrightarrow N$.

\subsection
{Periods}
The subspace $\Hom(\pi,\sut)^s$ maps to a (real) symplectic
submanifold $X_U\subset X$, and the corresponding real flows define
Hamiltonian systems which have been studied 
in~\cite{Inv}, and are
closely related to periodic flows studied by
Jeffrey-Weitsman~\cite{JeffreyWeitsman1,JeffreyWeitsman2,JeffreyWeitsman3,
JeffreyWeitsman4}.  The orbits of these flows are all closed, although
the period of the orbit varies with the value of $f$.

\begin{proposition}
Let $[\phi]\in X_U$ and let $\alpha\in\pi$ be represented by a simple loop $A$.
If $f_\alpha([\phi]\neq \pm 2$, the period of the 
trajectory of the flow of $\Ham(f_\alpha)$ at
$[\phi]$ equals
\begin{equation*}
\frac{4\pi}{\sqrt{4 - f_\alpha([\phi])^2}}
\end{equation*}
if $A$ is nonseparating and 
\begin{equation*}
\frac{2\pi}{\sqrt{4 - f_\alpha([\phi])^2}}
\end{equation*}
if $A$ separates.
If $f_\alpha([\phi])=\pm 2$, 
then $[\phi]$ is fixed under the flow of $\Ham(f_\alpha)$.
\end{proposition}

\begin{proof}
Let $P\in\sut$. Apply an inner automorphism of $\sut$ to assume that
$P$ is diagonal:
\begin{equation}\label{eq:defA}
P = \bmatrix e^{i\theta} & 0 \\  0 & e^{-i\theta}
\endbmatrix.
\end{equation}
Then 
\begin{equation*}
f(P) = 2 \cos(\theta), 
\qquad F(P) =   \bmatrix i\sin(\theta) & 0 \\  0 & -i\sin(\theta) \endbmatrix
\end{equation*}
with corresponding one-parameter subgroup 
\begin{equation*}
\zeta_t = \exp(tF(P)) =  \bmatrix e^{i\sin(\theta)t} & 0 \\  
0 & e^{-i\sin(\theta)t}  \endbmatrix
\end{equation*}
Since the $\pslt$-action on $\hpg^s$ is proper and free,
the quotient map is a principal $\pslt$-fibration:
\begin{equation*}
\pslt \longrightarrow \hpg^s  \longrightarrow X.
\end{equation*}
Thus the trajectory of the flow $\xi^{\alpha}_t$ 
projects diffeomorphically to the trajectory of the flow of 
$\Ham(f_\alpha)$ on $X$.  
The period equals the infimum of all $t>0$ such that 
\begin{equation*}
\xi^{\alpha}_t (\phi) = \phi. 
\end{equation*}

If $\alpha$ is nonseparating, then \eqref{eq:nonsep} implies that $T$ equals
the infimum of all $t>0$ such that 
\begin{equation*}
\zeta_t\phi(\beta) = \phi(\beta).  
\end{equation*}
By \eqref{eq:zetat},
\begin{equation*}
T = \frac{2\pi}{\sin(\theta)} = 
\frac{4\pi}{\sqrt{4 - f_\alpha([\phi])^2}}
\end{equation*}
If $\alpha$ separates, then \eqref{eq:sep} implies that $T$ equals
the infimum of all $t>0$ such that 
\begin{equation*}
\zeta_t\phi(\beta)\zeta_{-t} = \phi(\beta).  
\end{equation*}
for all $\beta\in\pi_1(M_-)$. 
If 
\begin{equation*}
\phi(\beta) \in \{\pm \Id\} = \mathsf{center}(\sut) 
\end{equation*}
for all $\beta\in\pi_1(M_-)$, then $\phi(\alpha) = 1$.
Since $\phi\in\hpg^s,$ there exists $\beta\in\pi$
such that $\phi(\beta)\neq 1$.
The period $T$ is the infimum of all $t>0$ such that
$\zeta_t = \pm \Id$.
By \eqref{eq:zetat},
\begin{equation*}
T = \frac{\pi}{\sin(\theta)} = 
\frac{2\pi}{\sqrt{4 - f_\alpha([\phi])^2}}
\end{equation*}
\end{proof}

\subsection
{Dehn twists}\label{sec:Dehntwists}
Closely related to the twist flows are the maps of $X_U$ induced by
Dehn twists. Recall that the {\em Dehn twist\/} about a simple loop
$A\subset M$ is the diffeomorphism supported in a tubular neighborhood
of $A$. If 
\begin{align*}
\psi: [0,1]\times S^1 & \hookrightarrow M \\
      \{\frac12\} \times S^1 & \hookrightarrow A 
\end{align*}
is such a tubular neighborhood, then the Dehn twist $\tau_A$ about $A$ is 
$\psi$-related to the diffeomorphism of $[0,1]\times S^1$
(restricting to the identity on the boundary) defined by:
\begin{equation*}
(s,e^{ i\theta}) \longmapsto (s, e^{2\pi i s}e^{i\theta}).
\end{equation*}
The action of $\tau_A$ on the fundamental group is given by:
\begin{equation}\label{eq:nonsepDehn}
(\tau_A)_* \ (\phi): \gamma \longmapsto 
\begin{cases} 
\phi(\gamma) & \text{~if~$\gamma\in\pi_1(M|A)$} \notag \\
\phi(\beta )\ \phi(\alpha) & \text{~if~$\gamma=\beta$.}\end{cases}
\end{equation}
if $A$ is nonseparating, and by:

\begin{equation}\label{eq:sepDehn}
(\tau_A)_* \ (\phi): \gamma \longmapsto 
\begin{cases}
\phi(\gamma) 
& \text{~if~$\gamma\in\pi_1(M_+)$} \notag\\
\phi(\alpha)
\phi(\gamma)
\phi(\alpha)^{-1}
& \text{~if~$\gamma\in\pi_1(M_-)$} \end{cases}
\end{equation}
if $A$ separates.

Let $P$ be as in \eqref{eq:defA}.
Comparing \eqref{eq:defA} with \eqref{eq:zetat}, $\zeta_t = P$ for
\begin{equation*}
\theta = \sin(\theta) t, 
\end{equation*}
that is,
\begin{equation}\label{eq:timeone}
t = \frac{\theta}{\sin(\theta)} = \frac{\cos^{-1}\big(f(P)/2\big)}
{2\pi\sqrt{4 - f(P)^2}}.
\end{equation}
Combining \eqref{eq:sepDehn}, \eqref{eq:sepDehn} with 
\eqref{eq:timeone} implies:

\begin{proposition}
Let $\alpha\in\pi$ correspond to a simple closed curve $A\subset M$.
The map $(\tau_A)_*$ on $\hpg$ induced by Dehn twist about $A$ 
equals the time-$t$ map of the flow generated by $\Ham(f_\alpha)$, where 
\begin{equation*}
t = \frac{2\cos^{-1}(f_\alpha([\phi])/2)}{\sqrt{4 - f_\alpha([\phi])^2}}.
\end{equation*}
\end{proposition}
In terms of a parametrization $\R/\Z\longrightarrow X_U$ of this
trajectory, $(\tau_A)_*$ acts by translation of
\begin{equation*}
\frac{\theta}{2\pi} = 
\frac{\cos^{-1}\big(f(P)/2\big)}{2\pi}
\end{equation*}
which has infinite order for almost every value of $f(A)$. 

By reparametrizing the flow, we obtain Hamiltonian flows whose
time-one map is the identity map or the Dehn twist. 
Jeffrey-Weitsman~\cite{JeffreyWeitsman2} (\S 5.1),
consider flows of Hamiltonians given by invariant functions
\begin{equation*}
\theta(P) = \cos^{-1}\bigg(\frac{f(P)}2\bigg)
\end{equation*}
which define $S^1$-actions, but are undefined at $P=\pm \Id$. 
Thus their flows are only defined on the 
on the dense open subset where $\phi(P)$ is not central. 
On the other hand their flows are periodic with period $2\pi$
if $A$ is nonseparating, and period $\pi$ if $A$ separates.
(Actually they work with the symplectic form which is $1/(4\pi^2)$
of ours, to obtain a 2-form with integral cohomology class.
Thus their periods are 
$1/(2\pi)$ 
if $A$ is nonseparating and 
$1/(4\pi)$ if $A$ separates.)

\subsubsection
{Orbits of complex twist flows}
On the $\slt$-character variety $X$, the Hamiltonian vector field
$\Ham(f_\alpha)$ generates a holomorphic $\C$-action. \eqref{eq:unipC}
implies that if $f(A) = \pm 2$ and $A\neq\pm \Id$ (that is, $A$ is
$\pm$-unipotent), then the trajectory defines an injective holomorphic
map 
\begin{align}
\C & \longrightarrow X \label{eq:unipCtraj}\\
t  & \longmapsto \xi^{\alpha}_t([\phi]). \notag
\end{align}
If $f(A)\neq\pm 2$, then \eqref{eq:diagC}
implies that $\zeta_t=1$ whenever $t\in \ t_0\Z$ where
\begin{equation}\label{eq:defnt0}
t_0 := \frac{4\pi i}{\lambda - \lambda^{-1}} =
\frac{4\pi i}{\big(f(A)^2 - 4\big)^{1/2}}
\end{equation}
(which is well-defined only up to sign).
In this case the trajectory of $[\phi]$ is the image of the holomorphic
embedding 
\begin{align}
\C^* &\longrightarrow X \label{eq:cxtraj} \\
e^z   &\longmapsto \xi^{\alpha}_{\big(2z\;/\;(\lambda -\lambda^{-1})\big)}\; 
\big([\phi]\big).\notag
\end{align}

\subsubsection
{Dehn twists}
At a point $[\phi]\in X$ where $f_\alpha(\phi) = \pm 2$, then either
\begin{itemize}
\item $\phi(\alpha) = \pm \Id$ is a central element and 
$[\phi]$ is fixed under the entire
$\C$-action, or 
\item $\phi(\alpha)$ is $\pm$ a parabolic element and \eqref{eq:unipC} implies
$\phi(\alpha) = \pm \exp F(\alpha)$.
\end{itemize}
In the latter case the embedding \eqref{eq:unipCtraj} is equivariant
with respect to the action of $\Z$ by translation on $\C$ and
the $\Z$-action generated by $(\tau_A)_*$ on $X$.

Suppose that $f_\alpha([\phi])\neq \pm 2$.  In that case $\zeta_t = P$ 
precisely when $t \equiv t_1\big(\operatorname{mod}t_0\big)$, where
\begin{equation*}
t_1 =  \frac{\log(\lambda)}{\lambda - \lambda^{-1}} =
\frac
{\log\big( (f(A)\pm(f(A)^2-4)^{1/2})/2\big)} 
{\big(f(A)^2-2\big)^{1/2}}.
\end{equation*}
(The two choices for $(f(A)^2-4)^{1/2}$ differ by sign
but determine equal values for $t_1$.)
In that case the embedding \eqref{eq:cxtraj} is equivariant
with respect to the actions generated by multiplication by
\begin{equation*}
\lambda = \frac{(f(A)\pm(f(A)^2-4)^{1/2})}2
\end{equation*}
on $\C^*$ and by $(\tau_A)_*$ on $X$.

\subsubsection
{Fenchel-Nielsen twist flows for $G=\sltr$}
Since the complex trace form $\B$ on $\gslt$ restricts to the trace form
on $\gsltr$, the complex-symplectic structure $\Omega$ on $X$ 
restricts to the symplectic structure on $\Teich$ defined by the
trace form of $\sltr$. By \cite{Nature} this symplectic structure
equals (-2) the Weil-Petersson K\"ahler form.

Inside the complex twist flows are the 
{\em Fenchel-Nielsen twist flows\/} when the holonomy is hyperbolic.
Suppose that $\phi\in\Hom(\pi,\sltr)$ is a discrete embedding
(that is, a {\em Fuchsian representation}). The subset of $X$ consisting
of equivalence classes of Fuchsian representations corresponds
bijectively to the space of marked hyperbolic structures on $M$,
that is, the {\em Teichm\"uller space $\Teich$.\/} In that case, for every
$\Id\neq \alpha\in\pi$, the element $\phi(\alpha)\in\sltr$ is hyperbolic.
Geometrically it is a {\em transvection\/} along a geodesic $\gamma$, and
is conjugate to a diagonal matrix
\begin{equation*}
P = \pm \bmatrix e^{l/2} & 0 \\ 0 & e^{-l/2} \endbmatrix.
\end{equation*}
$l$ is the distance $\phi(\alpha)$ moves points along $\gamma$.
In the hyperbolic structure on $M$, the element $\alpha$ corresponds
to a homotopy class of closed loops; $\gamma$ corresponds to the 
unique closed geodesic $A$ freely homotopic to a free loop in $\alpha$ and
$l$ is the length of this geodesic.

The character function is related to the length function by:
\begin{equation*}
f = 2\cosh(l/2).  
\end{equation*}
Thus the invariant function $l$ 
(which is only defined on the subset of $\sltr$ consisting of hyperbolic 
elements) defines functions 
\begin{equation*}
l_\alpha: \Teich \longrightarrow \R_+. 
\end{equation*}
As in \cite{Inv}, the corresponding variation is
\begin{equation*}
L(P) = \pm \bmatrix 1 & 0 \\ 0 & -1  \endbmatrix.
\end{equation*}
with corresponding one-parameter subgroup
\begin{equation*}
\zeta_t =  \bmatrix e^{t} & 0 \\ 0 & e^{-t} \endbmatrix,
\end{equation*}
the one-parameter subgroup of transvections moving at speed $2$ along the 
geodesic. 

The inner product of these infinitesimal transvections can be computed
as follows. The matrices
\begin{equation*}
L_1 = \bmatrix 0 & 1 \\ 1 & 0  \endbmatrix,
L_2 = \bmatrix -\sin(\psi) & \cos(\psi) 
\\ \cos(\psi) & \sin(\psi)  \endbmatrix 
\end{equation*}
represent infinitesimal transvections about axes which intersect with
angle $\psi$ and
\begin{equation*}
\B(L_1, L_2) = \tr(L_1L_2) = 2 \cos(\psi). 
\end{equation*}
Thus (as in \S 3.9 of \cite{Inv}), 
Wolpert's cosine formula~\ref{thm:WolpertDeriv}, follows from
Proposition~\ref{prop:PoissonBracket}. 

If $\alpha$ corresponds to a simple closed curve $A$, then the
corresponding twist flow $\eta^{\alpha}_t$ corresponds 
to the {\em  Fenchel-Nielsen twist flow\/} on 
$\Teich$ defined geometrically as 
follows. (Compare Abikoff~\cite{Abikoff}, Buser~\cite{Buser} and 
Harvey~\cite{Harvey}.) 
Represent $A$ by the (necessarily simple) closed
geodesic. Split $M$ along $A$, and identify the boundary components
$A_+,A_-$ of $M|A$ by a translation of length $t$ in the positively
oriented direction. This is a well-defined path of marked hyperbolic
structures in $\Teich$; see Wolpert~\cite{FenchelNielsen} for
details. In \cite{FenchelNielsen}, Wolpert proves that the Fenchel-Nielsen
twist flow is Hamiltonian with respect to the Weil-Petersson K\"ahler
form on $\Teich$, with Hamiltonian potential the length function
$l_\alpha$~\cite{Symplectic}.

The Fenchel-Nielsen twist flow $\eta^{\alpha}$ and the twist flow
$\xi^{\alpha}$ are re\-paramet\-ri\-za\-tions of one another. Specifically,
\begin{equation*}
\xi^{\alpha}_t =  \eta^{\alpha}_{\sinh(l_\alpha)}.
\end{equation*}

As in \cite{Nature} (see also \cite{Inv}), the Weil-Petersson K\"ahler
form equals $-1/2$ the symplectic structure defined by the trace
form $\B$ restricted to $\sltr$, implying the Fenchel-Nielsen flow is
Hamiltonian. Since the Fenchel-Nielsen twist flow is Hamiltonian for
the length function $l_\alpha$, the derivative of the length function
along the twist vector field equals the
Poisson bracket $\{l_\alpha,l_\beta\}$, which can be computed
by  Proposition~\ref{prop:PoissonBracket}. 
We thus obtain Wolpert's Derivative Formula:

\begin{thm}[(Wolpert~\cite{Symplectic,FenchelNielsen}]
\label{thm:WolpertDeriv}
Let $\alpha,\beta\in\pi$ where $\alpha$ is represented by a simple closed
curve. Then the derivative of the length function $l_\beta$ with respect
to the Fenchel-Nielsen twist flow equals the sum
\begin{equation*}
\sum_{i=1}^k\  \cos (\theta_i)
\end{equation*}
where $\alpha$ and $\beta$ are represented by closed geodesics,
$p_1,\dots,p_k$ are their intersection points, and $\theta_i$ is
the angle from $\alpha$ to $\beta$ at $p_i$.
\end{thm}

In particular the action of the Dehn twist $\tau_A$ equals the
time-$t$ map of the Fenchel-Nielsen twist flow at time $t = l/2$.

\section{Abelian Hamiltonian Actions}
The classical Fenchel-Nielsen coordinates on $\Teich$ can be
interpreted as a moment map arising from a {\em pants
decomposition,\/} that is, a maximal collection of disjoint simple
closed curves $\pp = \{\alpha_1,\dots,\alpha_N\}$ on $M$ which are each
homotopically nontrivial and mutually nonhomotopic. If $M$ has genus $g$,
then $N=3g-3$ and the length functions define a map
\begin{equation*}
l_\pp = (l_{\alpha_1}, \dots, l_{\alpha_N}):
\Teich \longrightarrow (\R_+)^N
\end{equation*}
which is the moment map for a Hamiltonian $\R^N$-action. This
$\R^N$-action is proper and free, and choosing a cross-section
(a left-inverse to $l_\pp$)
\begin{equation*}
\sigma:  (\R_+)^N \longrightarrow \Teich
\end{equation*}
defines a diffeomorphism
\begin{align}\label{eq:FNmap}
(\R_+)^N \times \R^N & \longrightarrow \Teich \\
(\lambda,t) & \longmapsto
\eta^{\alpha_1}_{t_1}\dots\eta^{\alpha_N}_{t_N}
\big(\sigma(\lambda)\big) \notag
\end{align}
The map $l_\pp$ defines the {\em action variables $\lambda$\/} while
the coordinates $t\in\R^N$ 
are the {\em angle variables.\/}
Indeed, Wolpert shows~\cite{WP} that this completely integrable system
defines a symplectomorphism with $(\R_+)^N \times \R^N$:
the Fenchel-Nielsen coordinates are Darboux coordinates with respect
to the symplectic structure 
$\omega_{WP}$ 
defined by the Weil-Petersson
K\"ahler form:
\begin{equation*}
\omega_{WP} = \sum_{i=1}^N d\lambda_i \wedge dt_i
\end{equation*}
This section extends this theory from $\sltr$ to $\slt$.

\subsection{Pants decompositions}

By Corollary~\ref{cor:disjoint}, the complex Hamiltonian vector fields
of $f_\alpha$ and $f_\beta$ Poisson commute if $\alpha$ and $\beta$
are represented by disjoint curves. Thus for a family 
of disjoint simple closed curves
\begin{equation*}
\pp = \{A_1,\dots,A_N\}
\end{equation*}
the corresponding complex twist flows of $\Ham{f_{\alpha_i}}$ generate
a $\C^n$-action.  Suppose these curves are each homotopically
nontrivial and mutually nonhomotopic.  The resulting map
\begin{equation*}
\tau_\pp: X \longrightarrow \C^N
\end{equation*}
is a moment map for a complex-Hamiltonian action of $\C^N$ on
the complex-symplectic manifold $X$.

Suppose that $\pp$ is {\em maximal,\/} that is $N=3g-3$, in which case
each component 
\begin{equation*}
P \subset M|(\cup_{i=1}^N A_i) 
\end{equation*}
is homeomorphic
to a 3-holed sphere. In that case Fricke's theorem~\cite{TopComps,Fricke}
implies that the generic inverse image
$(\tau_\pp)^{-1}(z)$ is a single $\C^N$-orbit. Specifically, let 
$(C^N)^s$ denote the subset consisting of 
$z\in\C^N$ such that,  for every $(i,j,k)$ for which
$\alpha_i,\alpha_j,\alpha_k$ bound a 3-holed  sphere $P_{ijk}$, 
\begin{equation*}
z_{\alpha_i}^2 + z_{\alpha_j}^2 +z_{\alpha_k}^2 -  
z_{\alpha_i} z_{\alpha_j} z_{\alpha_k} \neq 4.
\end{equation*}
This expresses the condition that the restriction  
$\phi|_{\pi_1(P_{ijk})}$ is irreducible; Fricke's theorem asserts
the triple
\begin{equation*}
\big(f(\phi(\alpha_i), f(\phi(\alpha_j), f(\phi(\alpha_k)\big)\in\C^3 
\end{equation*}
determines the equivalence class of such an irreducible representation.
If $z\in (\C^N)^s$, the fiber $(\tau_\pp)^{-1}(z)$ is a single $\C^N$-orbit. 

The subgroup $\gp$ of $\Gamma$ preserving the pants decomposition
is the free abelian group generated by the Dehn twists $\tau_{\alpha_i}$
for $i=1,\dots,N$. The resulting $\Z^n$-action lies in the $\C^N$-action
by the formulas in \S\ref{sec:Dehntwists}. Namely, at a point
$[\phi]\in X$ the orbit is an embedded product
\begin{equation*}
\prod_{i=1}^N G_i
\end{equation*}
where 
\begin{equation*}
G_i = \begin{cases} 
\{1\} &~\text{if $\alpha_i([\phi]) = \pm \Id$} \\ 
\C &~\text{if $f_{\alpha_i}([\phi]) = \pm 2$ and $\alpha_i([\phi]) \neq \pm \Id$} \\ 
\C^* &~\text{if $f_{\alpha_i}([\phi]) \neq \pm 2$} \end{cases}.
\end{equation*}
For each $i$ with with $f_{\alpha_i}(\phi)\neq 2$, let
\begin{equation*}
\lambda_i = 
\frac{f_{\alpha_i}(\phi) 
\pm(f_{\alpha_i}(\phi) ^2-4)^{1/2})}2.
\end{equation*}
Then $(k_1,\dots,k_N)\in\Z^N$ acts by translation by $k_i$ on the $i$-th
factor $G_i\approx\C$ for each $i$ with $f_{\alpha_i}(\phi) = \pm 2$
and by multiplication by $\lambda_i$ on each factor
$G_i\approx\C^*$ with $f_{\alpha_i}(\phi) \neq \pm 2$.

\subsection{The moment map as an ergodic decomposition and
its holomorphic extension}
The action of $\gp$ on $X_U$ is discussed in 
~\cite{Erg}. In this case the restriction
\begin{equation*}
f_\pp : X_U \longrightarrow [-2,2]^N 
\end{equation*}
is a moment map for a Hamiltonian $\R^n$-action whose orbits
are all closed. 
These orbits are tori, invariant under the $\gp\cong\Z^N$-action.
Almost every fiber has dimension $N$, and 
the action is equivalent to an $\Z^N$-action on $T^N$ by translation.
For almost all level sets, this action is ergodic. It follows that
the moment map $f_\pp$ for the Hamiltonian $\R^N$-action is also
the {\em ergodic decomposition\/} for the $\Z^N$-action (\cite{Erg},
Theorem~2.2):

\begin{proposition}\label{prop:ergdecom}
Let $h:X_U\longrightarrow\R$ be a 
$\gp$-invariant measurable function.
Then there exists a measurable function 
$\psi:[-2,2]^N\longrightarrow\R$
such that $h = \psi\circ f_\pp$ almost everywhere.
\end{proposition}

\begin{thm} Let $h:X\longrightarrow\cpo$ be a meromorphic function
which is invariant under $\gp$. Then there exists a meromorphic
function $H:\C^\pp\longrightarrow\cpo$ such that $h = H\circ\tau_\pp$.
\end{thm}
\begin{proof}
Let $G$ denote the complex-Hamiltonian $\C^\pp$-action on $X$
and let $G_U$ denote the Hamiltonian $\R^\pp$-action on $X_U$.
By Proposition~\ref{prop:ergdecom}, the measurable function $h$ must
be almost everywhere constant on the $G_U$-orbits on $X_U$. Thus
for each $t\in \tau_\pp(X_U)$, the function $h$ is constant on 
     the preimage $\tau_\pp^{-1}(t)\cap X_U$ and:
\begin{itemize}
\item 
The preimages of the restriction of $\tau_\pp$ to $X_U$ are the
$G_U$-orbits in $X_U$;
\item 
Each $G_U$-orbit in  $\tau_\pp^{-1}(t)$  is $\C$-Zariski dense in its 
$G$-orbit;
\item 
$G$ acts transitively on $\tau_\pp^{-1}(t)$.
\end{itemize}
Hence $h$ is constant on each $\tau_\pp^{-1}(t)$. for 
$t\in \tau_\pp(X_U)$.

Now we find a holomorphic section of $\tau_\pp$.
There exists a connected open neighborhood $W$ of $X_U$, a Zariski-closed
subset $Z\subset X$ of positive codimension, and a holomorphic map
\begin{equation*}
W-Z \stackrel{\sigma}\longrightarrow \tau_pp^{-1}(W-Z) 
\end{equation*}
so that the composition
\begin{equation*}
W-Z \stackrel{\sigma}\longrightarrow \tau_pp^{-1}(W-Z)  
\stackrel{\tau_\pp}\longrightarrow W-Z
\end{equation*}
equals the identity. Then the meromorphic function
\begin{equation*}
h - h\circ\sigma\circ\tau_\pp 
\end{equation*}
on $\tau_pp^{-1}(W-Z)$ vanishes on each fiber
$\tau_pp^{-1}(t)$ for $t\in W- Z \cap X_U$.
Since $W- Z \cap X_U$ is totally real, it follows that 
$h - h\circ\sigma\circ\tau_\pp$ vanishes on all of $W-Z$.
Since $W-Z$ is dense in $W$, this function vanishes on $W$.
Since $W$ is nonempty and open, this function vanishes on all of $X$.
Thus $h$ is constant on the fibers of $\tau_\pp$ as desired.
\end{proof}

\section{Deformation spaces of $\cpo$-structures}
A {\em $\cpo$-structure\/} on $M$ is a geometric structure with local
coordinate charts mapping to $\cpo$ with coordinate changes in
the group $\pslt$. The holonomy mapping $\hol$ maps the deformation space
$\cpo(M)$ locally biholomorphically to the character variety $X-X_U$.
Therefore $\Omega$ induces a complex-symplectic structure $\hol^*\Omega$
on $\cpo(M)$. On the other hand, $\cpo(M)$ is classically the total space
of an affine
bundle over $\Teich$ whose associated vector bundle is the cotangent bundle
$T^*\Teich$. The canonical exact complex-symplectic structure
on $T^*\Teich$ defines a complex-symplectic structure on $\cpo(M)$.
Kawai~\cite{Kawai} has proved this complex-symplectic structure
equals $\hol^*\Omega$:

\begin{thm*} {\rm (Kawai)} 
The complex-symplectic structure on $\cpo(M)$ induced
from the holomorphic cotangent bundle structure on 
     $T^*\Teich$ equals the
complex-symplectic structure induced from the complex-symplectic structure
$\Omega$ on $X$ by 
\begin{equation*}
\cpo(M) \stackrel{\hol}\longrightarrow \hpgg.
\end{equation*}
\end{thm*} 

The complex twist flows on $X$ lift to complex twist flows on
$\cpo(M)$. We relate the complex twist geometry on the open subset
$\qf\subset\cpo(M)$ comprising {\em quasi-Fuchsian structures\/} to
the quakebends studied by Epstein-Marden~\cite{EpsteinMarden},
Kourouniotis~\cite{Kourouniotis1,Kourouniotis2,Kourouniotis3},
McMullen~\cite{McMullen}, Platis~\cite{Platis}, Series~\cite{Series},
and Tanigawa~\cite{Tanigawa}.

\subsection{$\cpo$-structures}
Let $M$ be a smooth surface. A {\em $\cpo$-atlas\/} on $M$
consists of 
\begin{itemize}
\item An open covering $\uu$ of $M$ by {\em coordinate patches\/};
\item For each coordinate patch $U\in\uu$, a {\em coordinate chart\/}
\begin{equation*}
\psi_U:U \longrightarrow \cpo
\end{equation*}
which is a diffeomorphism of $U$ onto its image $\psi_U(U)$;
\item For each connected open subset $C\subset U\cap V$ of the intersection 
of coordinate patches $U,V\in\uu$, a transformation $g_C\in\pslt$
such that
\begin{equation*}
g_C \circ \psi_U|_{U\cap V} = \psi_V|_{U\cap V}.
\end{equation*}
The {\em coordinate change\/} $g_C$ is unique.
\end{itemize}
Such atlases are partially ordered by inclusion. A {\em $\cpo$-structure\/}
on $M$ is a maximal $\cpo$-atlas. A {\em $\cpo$-manifold\/} is a manifold
with a $\cpo$-structure.

Suppose $M$ and $M'$ are $\cpo$-manifolds. A mapping
$\phi:M\longrightarrow M'$ is a {\em $\cpo$-mapping\/} if and only if
for each coordinate chart $(U,\psi_U)$ in $M$ and each coordinate
chart $(U',\psi'_{U'})$ in $M'$ the composition
\begin{equation*}
(\psi'_{U'})\circ \phi\circ(\psi_U)^{-1}: \;
\psi_U\big(U \cap \phi^{-1}(U')\big) \;\longrightarrow\;
\psi'_{U'}\big(\phi(U)\cap U'\big)
\end{equation*}
extends to an element of $\pslt$ on each component of 
$\psi_U(U\cap\phi^{-1}(U')$. 
A $\cpo$-mapping is necessarily a local diffeomorphism.

If $f:M\longrightarrow M'$ is a local diffeomorphism, then every
$\cpo$-structure on $M'$ determines a unique $\cpo$-structure on $M$
such that $f$ is a $\cpo$-mapping. In particular every covering space
of a $\cpo$-manifold is a $\cpo$-manifold, and its group of covering
transformations is realized by a group of $\cpo$-automorphisms.

If $M$ is a simply-connected $\cpo$-manifold, then any chart on $M$
extends to a globally defined {\em developing map\/} 
\begin{equation*}
\dev:M \longrightarrow \cpo 
\end{equation*}
which is a $\cpo$-mapping. This mapping is unique up to left-composition
with transformations in $\pslt$. 

We denote by $\Sigma$ the Riemann surface whose underlying manifold is $M$
for which $\dev$ is holomorphic.

In general the universal covering space $\tM$ of a $\cpo$-manifold
admits a developing map into $\cpo$. The group $\pi_1(M)$ of covering
transformations of $\tM\longrightarrow M$ admits a representation
$\rho:\pi_1(M)\longrightarrow\pslt$ such that the diagram
\begin{equation*}
\begin{CD}
\tM @>\dev>> \cpo \\
@V{\gamma}VV    @VV{\rho(\gamma)}V \\
\tM @>\dev>> \cpo \end{CD}
\end{equation*}
commutes for every $\gamma\in\pi_1(M)$. The representation $\rho$ is
called the {\em holonomy representation.\/} The {\em developing
pair\/} $(\dev,\rho)$ is unique up to the action of $\pslt$ (by
left-composition with $\dev$ and conjugation on $\rho$).  The developing
map globalizes the coordinate charts and the holonomy representation
globalizes the coordinate changes.

\subsubsection{Deformation spaces}
Let $S$ be a fixed topogical surface with fundamental group $\pi=\pi(S)$. 
A {\em marked $\cpo$-structure \/} on $S$ consists of a $\cpo$-manifold
$M$ and a homeomorphism $f_M:S\longrightarrow M$. Two marked $\cpo$-structures
$(M,f_M)$ and $(M',f_{M'})$ are {\em equivalent\/} if and only if there
exists a $\cpo$-isomorphism $\phi:M\longrightarrow M'$ such that
\begin{equation*}
\phi\circ f_M \simeq f_{M'}.
\end{equation*}
Let $\cpo(S)$ denote the set of equivalence classes of marked
$\cpo$-structures on $S$. Holonomy defines a mapping
\begin{equation*}
\hol:  \cpo(S) \longrightarrow \hpgg.
\end{equation*}
As in \cite{Geost} (see also Earle~\cite{Earle},
Hubbard~\cite{Hubbard} and Kapovich~\cite{Kapovich}), 
there is a natural topology on 
$\cpo(S)$ such that $\hol$ is a local homeomorphism. 
Furthermore Gallo-Kapovich-Marden~\cite{GKM} proved that
the image of $\hol$ is $X-X_U-X_D$, where $X_D\subset X-X_U$ is the
closed subset consisting of
equivalence classes $[\phi]$ of representations leaving invariant
a pair of points of $\cpo$ (these correspond to irreducible representations
for which an index-two subgroup acts reducibly).

\subsection{Holonomy and the conformal structure}

Since $\pslt$ acts holomorphically on $\cpo$, a $\cpo$-atlas on $M$ is
a {\em holomorphic atlas,\/} that is, an atlas for a complex structure
on $M$. Thus underlying every $\cpo$-manifold is a Riemann surface.  Recording
the complex structure underlying a $\cpo$-structure is a map
\begin{equation*}
\cpo(S) \stackrel{\Pi}\longrightarrow \Teich. 
\end{equation*}
A {\em projective structure on a Riemann surface $\Sigma$\/}
is a $\cpo$-structure whose underlying complex structure is $\Sigma$.

\pagebreak

\begin{lemma}
$\Pi$ is holomorphic.
\end{lemma}

\begin{proof}
We recall the classical description of $\cpo(S)$ as an affine bundle
$\mathcal{Q}(M)$ over $\Teich$ whose associated vector bundle equals
the cotangent bundle $T^*(\Teich)$.
(See for example \S 2.3 of Earle~\cite{Earle}, Hubbard~\cite{Hubbard} or 
Gunning~\cite{Gunning}). 
Fix a Riemann surface $\Sigma$
homeomorphic to $M$ and a projective structure $s_0$ on $\Sigma$. 
Given any 
other projective structure on $\Sigma$, its developing map is a
holomorphic map from the universal covering $\tilde\Sigma$ to
$\cpo$. Its Schwarzian derivative (in the local projective coordinates
defined by $s_0$) is a holomorphic quadratic differential which is
invariant under $\pi$, and therefore defines a holomorphic quadratic
differential on $\Sigma$. Furthermore this quadratic differential
completely determines the developing map up to an element of $\pslt$
of $\cpo$.  Thus the fiber $\Pi^{-1}([\Sigma])$ admits a simply
transitive action of the
complex vector space $H^0(\Sigma;\ooo(\kappa^2))$ 
comprising holomorphic quadratic
differentials on $\Sigma$.  The composition law for the Schwarzian
derivative implies that changing the origin $s_0$ effects a
translation of $Q(\Sigma)$, so $\Pi^{-1}(\Sigma)$ is an affine space
modelled on $H^0(\Sigma;\ooo(\kappa^2))$, which is the cotangent
space $T^*_{[\Sigma]}\Teich$ of $\Teich$ at $[\Sigma]$.

In particular $\cpo(S)$ inherits a complex structure, making it a holomorphic 
affine bundle over $\Teich$. In particular the projection $\Pi$ is
holomorphic in this complex structure.

By Lemma~2 of Earle~\cite{Earle} or Hubbard~\cite{Hubbard}, the
holonomy mapping $\hol$ is a local biholomorphism, and therefore this
complex structure on $\cpo$ induced by $\hol$ is the above complex
structure. Thus $\Pi$ is holomorphic with respect to the complex
structure induced by $\hol$.
\end{proof}
Kawai's theorem implies that $\Pi$ is a holomorphic Lagrangian fibration. 
As described in Hubbard~\cite{Hubbard}, 
the differential $d\Pi$ identifies with the map 
\begin{equation*}
T_{s_0}\cpo(S) \cong H^1\big(\Sigma;\gslt\big) \longrightarrow  
H^1\big(\Sigma;\ooo(\kappa^{-1})\big) \cong T_{[\Sigma]}\Teich
\end{equation*}
induced by the short exact sequence of sheaves
\begin{equation}\label{eq:exactsheaf}
0\longrightarrow  \gslt \stackrel{i}\longrightarrow \ooo(\kappa^{-1})
\stackrel{D_3}\longrightarrow \ooo(\kappa^{2}) 
\longrightarrow 0.
\end{equation}
Here $D_3$ is the map of holomorphic sheaves given in local holomorphic
coordinates by
\begin{equation*}
h(z) \,\frac{\partial}{\partial z} 
\longmapsto h'{}'{}'{}(z)\, dz^2.
\end{equation*}
Its kernel $\gslt$ is the locally constant sheaf of locally projective
vector fields (whose coefficients in affine coordinates are quadratic).
The monomorphism $i$ regards a locally projective vector field as a 
holomorphic vector field. Then
\eqref{eq:exactsheaf} induces the exact sequence
\begin{equation*}
0 \longrightarrow H^0\big(\Sigma,\ooo(\kappa^2)\big) 
\overset{\delta}\longrightarrow
H^1\big(\Sigma;\gslt\big) \overset{i_*}\longrightarrow  
H^1\big(\Sigma;\ooo(\kappa^{-1})\big) \longrightarrow 0. 
\end{equation*}
If $\alpha\in H^0\big(\Sigma,\ooo(\kappa^2)\big)$ is a holomorphic
quadratic differential  
and $\beta\in H^1\big(\Sigma;\gslt\big)$, then
\begin{equation}\label{eq:reciprocity}
\langle \delta(\alpha), \beta \rangle = 
\alpha\cdot i_*(\beta)
\end{equation}
where the first pairing $\langle,\rangle$ is the pairing
on $H^1\big(\Sigma;\gslt\big)$ induced by the complex
Killing form $\B$ on $\gslt$ and the second pairing is the
Serre duality pairing
\begin{equation*}
H^0\big(\Sigma,\ooo(\kappa^2)\big) \times
H^1\big(\Sigma,\ooo(\kappa^{-1})\big) \longrightarrow
H^1\big(\Sigma,\ooo(\kappa)\big) \cong \C.
\end{equation*}
This discussion implies that $\Pi$ is a Lagrangian fibration.
Suppose 
\begin{equation*}
\beta,\gamma\in 
\ker(d\Pi) =  \delta( H^0\big(\Sigma,\ooo(\kappa^2)\big) =
\ker(i_*)
\end{equation*}
Then $\gamma = \delta(\alpha)$ for some
$\alpha\in H^0\big(\Sigma,\ooo(\kappa^2)\big)$, and 
\eqref{eq:reciprocity} implies 
$\langle \gamma,\beta\rangle = 0$ 
since $i_*\beta = 0$.
Therefore $\ker(d\Pi)$ is isotropic. Since its dimension is half of
that of $H^1\big(\Sigma;\gslt\big)$, the fibers of $\Pi$ are Lagrangian.
For further details, see Kawai~\cite{Kawai}.

\subsubsection{Meromorphic functions on $\cpo(S)$}
\begin{proposition}\label{cor:nonconstant}
There exist nonconstant $\Gamma$-invariant meromorphic functions
on $\cpo(S)$.
\end{proposition}
\begin{proof}
The
{\em Riemann moduli space of curves\/} is the 
quotient $\Teich/\Gamma$, which by Knudsen~\cite{Knudsen}, is
a quasiprojective variety. Thus 
$\Teich/\Gamma$ admits nonconstant meromorphic functions and
$\Teich$ admits nonconstant $\Gamma$-invariant meromorphic functions. 

Let $\psi$ be a nonconstant $\Gamma$-invariant meromorphic function on 
$\Teich$. Composing with the projection 
\begin{equation*}
\Pi:\cpo(S) \longrightarrow \Teich
\end{equation*}
provides a nonconstant $\Gamma$-invariant meromorphic function
$\psi\circ\Pi$ on $\cpo(S)$.
\end{proof}

\subsection{Quasi-Fuchsian space}\label{sec:quasifuchsian}
A $\cpo$-structure is {\em Fuchsian\/} if and only if a 
developing map embeds $\tM$ in a hyperbolic plane
$\cho\subset\cpo$. Necessarily the holonomy representation is
conjugate to a {\em Fuchsian representation,\/} that is, a discrete
embedding $\pi\hookrightarrow\puoo$ (where $\puoo$ is the stabilizer
in $\pslt$ of the Poincar\'e disc $\cho$). However by a construction
of Maskit and Hejhal (see also Goldman~\cite{Fuchsian}) there exist
$\cpo$-structures with Fuchsian holonomy with surjective developing
maps. A {\em quasi-Fuchsian $\cpo$-structure\/} is a $\cpo$-structure
topologically conjugate to a Fuchsian structure, that is, there exists
a homeomorphism $h$ of $\cpo$ such that
\begin{equation*}
(h\circ\dev, h\circ\rho\circ h^{-1}) 
\end{equation*}
is the developing pair for a Fuchsian structure.
The limit set $\Lambda$ of $\rho(\pi)$ is the Jordan curve
$h(\partial\cho)$. The developing image $\dev(\tM)$ is
one component of the complement $\cpo-\Lambda$.

{\em Quasi-Fuchsian space\/} $\qf$ is the subset of $\cpo(S)$ corresponding
to quasi-Fuchsian structures.   The holonomy representation
$\rho$ is a {\em quasi-Fuchsian representation.\/} The
holonomy mapping
\begin{equation*}
\hol:\qf \longrightarrow X 
\end{equation*}
embeds $\qf$ as the open subset of $X$ consisting of conjugacy classes
of quasi-Fuchsian representations.  Quasi-Fuchsian
space $\qf$ is both the open subset of $\cpo(S)$ consisting of
quasi-Fuchsian structures and the open subset of $X-X_U$ consisting of
characters of quasi-Fuchsian representations.

The conformal structures of the quotients of $\cpo-\Lambda$ 
by $\phi(\pi)$ determine an ordered pair 
$\Teich\times\bTeich$; the first parameter is the marked
conformal structure underlying the $\cpo$-structure.
The {\em simultaneous uniformization theorem\/} of Bers~\cite{Bers}
asserts that the corresponding map
\begin{equation*}
\qf \stackrel{\approx}\longrightarrow\Teich\times\bTeich  
\end{equation*}
is a biholomorphism. It is evidently $\Gamma$-equivariant.

Since $\qf\hookrightarrow\cpo(S)$, Corollary~\ref{cor:nonconstant}
implies that nonconstant $\Gamma$-invariant meromorphic functions exist 
on $\qf$.

\subsection{Twist flows of $\cpo$-structures and grafting}
Let $\alpha\in\pi$. Then the composition
\begin{equation*}
\cpo(S) \stackrel{\hol}\longrightarrow X \stackrel{f_\alpha}\longrightarrow \C
\end{equation*}
defines a holomorphic function on $\cpo(S)$ and its complex-Hamiltonian
$\Ham(f_\alpha\circ\hol)$ is a holomorphic vector field on $\cpo(S)$.  If
$\alpha$ corresponds to a simple closed curve $A$, then there is a
{\em complex twist flow\/} 
$\tilde{\xi}^{\alpha}_t$ on $\cpo$ covering the
complex twist flow generated by the vector field $\Ham(f_\alpha)$ on $X$.

This twist flow is defined geometrically on Fuchsian structures as
follows (see \cite{Fuchsian},
Kourouniotis~\cite{Kourouniotis1,Kourouniotis2,Kourouniotis3}).  Let
$s_0$ be a Fuchsian $\cpo$-structure, so a developing map embeds $\tM$
as a geometric disc $\Delta\subset \cpo$. Then the simple closed curve
$A\subset M$ can be represented by a simple closed geodesic with
respect the Poincar\'e metric induced from that of $\Delta$.
Let $\ell(A)$ denote the Poincar\'e length of this geodesic.
A lift $\tilde A$ of $A$ to $\tM$ develops to a circular arc orthogonal
to $\partial\dev(\tM)$. The split surface $M|A$ inherits a $\cpo$-structure
whose boundary components develop to circular arcs.
To define $s_t := \tilde{\xi^{\alpha}}_t(s_0)$, insert an annulus $A_\theta$ 
with $\cpo$-structure into $M|A$.
The angular parameter $\theta$ is the imaginary part $\Im(t)$.
The annulus is a {\em $\theta$-annulus\/} 
in the sense of \S 2.12 of~\cite{Fuchsian}.
Choose a holomorphic universal covering map
\begin{equation*}
E:\C \longrightarrow \cpo-\Fix\big(\phi(\alpha)\big) 
\end{equation*}
which is periodic with period $2\pi i$.
(When $\Fix(\phi(\alpha)) = \{0,\infty\}$, then this is just
the exponential map $\exp:\C\longrightarrow \C^*$. For general $A$,
we compose $\exp$ with a projective transformation taking $\{0,\infty\}$
to $\Fix(\phi(\alpha))$.)

Define the {\em $\theta$-strip\/} 
\begin{equation*}
S_\theta := \R + i[0,\theta] \subset \C 
\end{equation*}
with $\cpo$-structure induced by $E$. The developing map $E$ for
$S_\theta$ is equivariant with respect to the $\Z$-action on $S_\theta$
generated by translation by $\ell(A)$ and the $\Z$-action on 
$\cpo-\Fix(\phi(\alpha))$. The {\em $\theta$-annulus\/} $A_\theta$ is
defined as the quotient $\cpo$-manifold of $S_\theta$. 

The {\em grafted manifold\/} $M(t)$ is obtained by
inserting a $\theta$-annulus into $M|A$. Choose one component of
$A_- \subset \partial(M|A)$ and attach $A_\theta$ to $A_-$ to obtain
a $\cpo$-manifold homeomorphic to $M|A$ with one boundary component
$A_+'$ corresponding to the component $A_+$ of $\partial(M|A)$ and another
boundary component $A'$ corresponding to the other component
of $\partial(A_\theta)$. 
Identify $A_+'$ to $A'$ by a transvection of displacement $\Re(t)$.

When $\theta=0$, this is just the Fenchel-Nielsen twist deformation,
obtained from $M|A$ by identifying the two components of $\partial(M|A)$
by a transvection of displacement $\Re(t)$.

Let $P\in\slt$. These holomorphic flows on $\hpg$ cover holomorphic
flows on $X$ which have been extensively studied on quasi-Fuchsian
space $\qf\subset X.$ On $\qf$ these flows geometrically correpond to
the {\em quakebends\/} or {\em complex earthquakes\/} discussed by
Epstein-Marden~\cite{EpsteinMarden},
Kourouniotis~\cite{Kourouniotis1,Kourouniotis2,Kourouniotis3},
McMullen~\cite{McMullen}, Platis~\cite{Platis},
Series~\cite{Series}, and Tanigawa~\cite{Tanigawa}. 
In particular we obtain the following result of Platis~\cite{Platis}:

\begin{thm}[(Platis~\cite{Platis}, Theorem 7, \S 2.2]
The complex twist vector field $\tilde\xi$ on $\qf$ is Hamiltonian with 
respect to the complex length function $l^{\C}$. 
\end{thm}

Just as the geodesic length function $l$ on the hyperbolic subset
of $\sltr$ and the angle function $\theta$ on $\sut - \{\pm\Id\}$
are more geometrically natural than the trace functions,
we consider the {\em complex length\/} ``function'' on $\slt$ defined by:
\begin{equation*}
2 \cosh\bigg(\frac{l^{\C}(P)}2\bigg) = f(P)
\end{equation*}
or equivalently
\begin{equation*}
l^{\C}(P) = 2 \log\bigg( 
\frac 
{ \big( f(P) \pm (f(P)^2-4)^{1/2} \big) }
{2} 
\bigg).
\end{equation*}
(This differs from Tan~\cite{Tan} whose complex length is half of ours.
Our definition is consistent with 
Kourouniotis~\cite{Kourouniotis1,Kourouniotis2,Kourouniotis3}, 
Platis~\cite{Platis} and Series~\cite{Series}.)
This function takes values in $(\C/4\pi i\Z)/\{\pm 1\}$, since
the logarithm is well-defined only modulo $4\pi i$ and the choice 
of square root introduces an ambiguity of sign.

Choose the branch $\tilde{l}_\C$ of the complex length
which is positive on the subset of hyperbolic elements in $\sltr$;
since $\qf$ is simply connected
\begin{equation*}
[\phi] \longmapsto 
\big(\tilde{l}^\C(\phi(\alpha_1)),\dots, \tilde{l}^\C(\phi(\alpha_N))\big)
\end{equation*}
defines a single-valued function
\begin{equation*}
l^{\C}|_\pp: \qf \longrightarrow  \C^N.
\end{equation*}
The restriction of $l^{\C}_\pp$ to $\Teich$ is the length function
$l_\pp:\Teich\longrightarrow (\R_+)^N$ defined by \eqref{eq:FNmap}.
The fibers of $l^{\C}|_\pp$ are orbits of the complex-Hamiltonian 
$\C^N$-action having $l^{\C}|_\pp$ as moment map. Tan~\cite{Tan} and 
Kourouniotis~\cite{Kourouniotis3} 
show that there exists a section $\sigma:\C^N\longrightarrow\qf$ 
and a neighborhood 
of $(\R_+)^N \times \{0\}$ in $\C^N\times\C^N$ 
which maps biholomorphically to $\qf$.

Korouniotis~\cite{Kourouniotis2} and Series~\cite{Series} derive formulas
for the derivative of the complex length functions along quakebend flows.
We briefly sketch how their formulas can be derived from 
Proposition~\ref{prop:PoissonBracket} as a Poisson bracket, referring
to \cite{Kourouniotis2,Series} for details.

Let $P\in\slt$ be loxodromic with invariant axis $a_P$.
Then the variation function $L^\C$ associated
to the invariant function $l^\C$ and the complex-orthogonal structure
$\B(X,Y) = \tr(XY)$  maps
a diagonal matrix to
\begin{equation*}
D_0 := \pm \bmatrix 1 & 0 \\ 0 & -1 \endbmatrix 
\end{equation*}
(which is not uniquely defined without restrictions to make the
function $l^\C$ single-valued).  By \eqref{eq:equivariance}, the
element of $\pslt$ corresponding to $L^\C(P)$ is the involution fixing
$a_P$. If $P_1,P_2\in\pslt$ are loxodromic, with complex distance
$z$ between their invariant axes, then a simple
calculation implies that $\B(L^\C(P_1),L^\C(P_2))$ equals $\cosh(z)$.
Applying Proposition~\ref{prop:PoissonBracket} to the complex length functions
$l^\C_\beta$ and $l^\C_\alpha$, where $\alpha$ corresponds to a simple
closed curve, we obtain the following extension of 
Theorem~\ref{thm:WolpertDeriv}:

\begin{thm}[(Korouniotis~\cite{Kourouniotis2}, Series~\cite{Series}]
Let $\alpha,\beta\in\pi$ where $\alpha$ is represented by a simple
closed curve $A$. Then the derivative of the complex length function
$l^C_\beta$ on $\qf$ with respect to the quakebend flow with respect
to $A$ equals the sum
\begin{equation*}
\sum_{i=1}^k\  \cosh (d_i)
\end{equation*}
where $\alpha$ and $\beta$ are represented by closed geodesics,
$p_1,\dots,p_k$ are their intersection points, and $d_i$ the complex
distance between the axes of $\phi_i(\alpha_i)$ to $\phi_i(\beta_i)$ at $p_i$,
where $\phi_i$ is a representative homomorphism from $\pi_1(M;p_i)$ and
$\alpha_i,\beta_i$ are elements of $\pi_1(M;p_i)$ representing $\alpha$ and
$\beta$.
\end{thm}
See Korouniotis~\cite{Kourouniotis2} and Series~\cite{Series} for details.

Tan~\cite{Tan} and 
Kourouniotis~\cite{Kourouniotis1,Kourouniotis2,Kourouniotis3} extend 
Fenchel-Nielsen coordinates to complex Fenchel-Nielsen coordinates on
$\qf$. (See also Series~\cite{Series} for another exposition.)

\end{document}